\documentclass[10pt]{article}
\usepackage{amssymb}
\usepackage{graphicx}
\usepackage{amsmath}
\usepackage{color}
\usepackage[a4paper]{geometry}

\newtheorem{thm}{Theorem}[section]
 \newtheorem{cor}[thm]{Corollary}
 \newtheorem{lem}[thm]{Lemma}
 \newtheorem{prop}[thm]{Proposition}
 
 \newtheorem{rem}[thm]{Remark}

 \numberwithin{equation}{section}

\def\hh{\!\!\!\!}

\def\Proof{\noindent\emph{Proof.} }
\def\qed{\hfill $\Box$ \smallskip}

\begin{document}
\title{\textbf{Direct and inverse problems for a third order self-adjoint differential operator with periodic boundary conditions and nonlocal potential}}
\author{Yixuan Liu$^1$, \qquad Jun Yan$^{2}$\thanks{Corresponding author.}}
\date{}
\maketitle
\begin{center}
$^1$ School of Science, Civil Aviation University of China, Tianjin, 300300, People's Republic of China, \\
$^2$ School of Mathematics, Tianjin University, Tianjin, 300354, People's Republic of China

E-mail: \texttt{yx\_liu@cauc.edu.cn, jun.yan@tju.edu.cn}

\today
\end{center}

\begin{abstract}
A third order self-adjoint differential operator with periodic boundary conditions and an one-dimensional perturbation has been considered. For this operator, we first show that the spectrum consists of simple eigenvalues and finitely many eigenvalues of multiplicity two. Then the expressions of eigenfunctions and resolvent are described. Finally, the inverse problems for recovering all the components of the one-dimensional perturbation are solved. In particular, we prove the Ambarzumyan-type theorem and show that the even or odd potential can be reconstructed by three spectra.


\end{abstract}

\maketitle

% This is the full title of the paper
% Use lowercase letters in title except for proper names
% Avoid equations in title if possible
% Do not use the \thanks{} command; use \extraline{} instead (see below).

%Insert `2000 Mathematics Subject Classification' numbers here:

\section{Introduction}
In this manuscript, we focus on a third order self-adjoint differential operator $L_{\alpha}=L_{\alpha}(v)$ in $L_{\mathbb{C}}^2(0,1)$ defined by
\begin{equation*}
(L_{\alpha}y)(x):=\mathrm{i}y^{\prime\prime\prime}(x)+\alpha\int_{0}^{1}y(t)\overline{v}(t)\mathrm{d}tv(x),\quad \alpha\in\mathbb{R},\quad v\in L_{\mathbb{c}}^{2}(0,1),
\end{equation*}
whose domain $\mathcal{D}(L_{\alpha})$ consists of the functions $y\in W_{3}^{2}(0,1)$ satisfying the periodic boundary conditions
\begin{equation}
y(1)=y(0),\quad y^{\prime}(1)=y^{\prime}(0),\quad y^{\prime\prime}(1)=y^{\prime\prime}(0).\label{bc}
\end{equation}
Here the nonlocal potential $(\alpha,v)$ is crucial in mathematical physics \cite{qst}, and the operator $L_{\alpha}$ can be considered as a one-dimensional perturbation of the operator $L_{0}y=\mathrm{i}y^{\prime\prime\prime}$.

The motivation for studying third order differential operator stems from its role for constructing the Lax L-A pair of nonlinear Boussinesq, Camassa-Holm and Degasperis-Procesi equations. These nonlinear evolution equations can be used to investigate the long waves in shallow water, nonlinear lattice waves and vibrations of a cubic string (see \cite{mckean1981cpam, cs, dpe} and the references therein). The direct and inverse problems for the third order self-adjoint differential operator on a finite interval have attracted considerable attention in recent years \cite{amour1999siam, amour2001siam, ugurlu, lsy, atheorem, VAZ2021}. In particular, Zolotarev \cite{VAZ2021} investigated the spectrum of operator $L_{\alpha}$ on $[0,l]$ with the boundary conditions $y(0)=y(l)=0$, $y^{\prime}(0)=y^{\prime}(l)$ and found that the nonlocal potential can be recovered by four spectra. In this paper, we generalize some results in \cite{VAZ2021} to the operator $L_{\alpha}$ with periodic boundary conditions (\ref{bc}). It should be noted that these generalizations are not common. Firstly, owing to the periodic boundary conditions (\ref{bc}), the calculation of the expression for the resolvent of the operator $L_0$ is more complicated (see Lemma \ref{rel0}). Secondly, there may exists triple eigenvalues due to the three coupled boundary conditions (\ref{bc}). In fact, a very crucial Lemma \ref{csd2npi} and a detailed computation in the proof of Lemma \ref{mef} show that the spectrum of $L_0$ consists of simple eigenvalues; then on the basis of these results and according to the expression for the resolvent of $L_\alpha$, we find the spectrum of $L_{\alpha}$ is composed of simple eigenvalues and finitely many eigenvalues of multiplicity two. Finally, we also mention that for the first time we prove the Ambarzumyan-type theorem for the operator $L_\alpha$ (see Theorem \ref{ath}) and show that the even potential ($\overline{v}(1-x)=v(x)$, see Theorem \ref{do}) or odd potential ($\overline{v}(1-x)=-v(x)$, see Corollary \ref{do1}) can be reconstructed by three spectra.

%For the operator $L_0$, one of the most interesting and important features of its spectrum is consist of simple eigenvalues $2n\pi$, $n\in\mathbb{Z}$, which simple property seems contrary to the coupled boundary conditions (\ref{bc}). Furthermore, its eigenfunctions $\mathrm{e}^{2n\pi\mathrm{i}x}$ form an orthnormal basis of $L_{\mathbb{C}}^2(0,1)$ and help us to deduce the expressions of resolvent and characterize the spectrum of the operator $L_{\alpha}$, which are essential for solving the following inverse problems.

This paper is structured as follows. Some basic notations and useful properties are given in section 2. In section 3, we prove that the spectrum of $L_{0}$ is simple and deduce the expressions of corresponding eigenfunctions and resolvent operator. The spectrum, eigenfunctions and resolvent operator of $L_{\alpha}$ are investigated in section 4. In section 5, we solve the inverse problems for reconstructing the nonlocal potential $(\alpha,v)$ under different circumstances.

%\begin{theorem}
%For $v\in L^{2}(0,1)$, $\Vert v\Vert=1$, $\alpha\in \mathbb{R}$, the non-local potential $\{\alpha, v\}$ for the operator $L_{\alpha}$ can be reconstructed by four spectra.
%If $\overline{v(1-x)}=\pm v(x)$, then the non-local potential $\{\alpha, v\}$ for the operator $L_{\alpha}$ can be reconstruct by three spectra.
%\end{theorem}
\section{Preliminaries}
%In this section, some useful notations and basic properties are introduced.
In the following, we use $\left<f,g\right>=\int_{0}^{1}f(x)\overline{g}(x)\mathrm{d}x$ to denote the inner product for $f,g\in L_{\mathbb{C}}^2(0,1)$. The norm for $f\in L_{\mathbb{C}}^2(0,1)$ is defined by $\Vert f\Vert=\left<f,f\right>^{\frac{1}{2}}$.
\begin{prop}
\label{omega}
Denote $\omega=\mathrm{e}^{\frac{2\pi}{3}\mathrm{i}}=-\frac{1}{2}+\frac{\sqrt 3}{2}\mathrm{i}$, then the following equalities hold:
\begin{eqnarray*}
&&\omega^2=-\frac{1}{2}-\frac{\sqrt 3}{2}\mathrm{i}=\overline{\omega},\quad \omega^{3}=1, \quad 1+\omega+\omega^{2}=0,\\\notag
&&1-\omega=\mathrm{i}\sqrt 3 \omega^{2},\quad\omega-\omega^{2}=\mathrm{i}\sqrt 3,\quad\omega^{2}-1=\mathrm{i}\sqrt 3 \omega.\label{omega}
\end{eqnarray*}
\end{prop}
Let
\begin{equation}
c(z)=\frac{1}{3}\sum_{k=1}^{3}\mathrm{e}^{\omega^{k}z},\quad
s(z)=\frac{1}{3}\sum_{k=1}^{3}\frac{1}{\omega^{k}}\mathrm{e}^{\omega^{k}z},\quad
d(z)=\frac{1}{3}\sum_{k=1}^{3}\omega^{k}\mathrm{e}^{\omega^{k}z} \label{csd}
\end{equation}
denote the fundamental solutions of equation $y^{\prime\prime\prime}=y$ determined by the initial conditions
\begin{equation*}
\left(
\begin{array}{ccc}
c(0) & s(0) & d(0) \\
c^{\prime }(0) & s^{\prime }(0) &d^{\prime }(0) \\
 c^{\prime \prime}(0) &s^{\prime \prime}(0) & d^{\prime \prime}(0)%
\end{array}%
\right) =I_{3}:=\left(
\begin{array}{ccc}
1 & 0 & 0 \\
0 & 1 & 0 \\
0 & 0 & 1%
\end{array}%
\right).
\end{equation*}
The following lemmas show some basic properties for the solutions $c(z)$, $s(z)$, $d(z)$.
\begin{lem}
\label{entire}
\rm{\cite[Lemma 1.1]{VAZ2021}} The functions $c(z)$, $s(z)$ and $d(z)$ in $(\ref{csd})$ are entire functions of exponential type, and satisfy the following relations:\\
\rm{(i)} $s^\prime(z)=c(z)$, $d^\prime(z)=s(z)$, $c^\prime(z)=d(z)$;\\
\rm{(ii)} $\overline{c(z)}=c(\overline{z})$, $\overline{s(z)}=s(\overline{z})$, $\overline{d(z)}=d(\overline{z})$;\\
\rm{(iii)} $c(\omega z)=c(z)$, $s(\omega z)=\omega s(z)$, $d(\omega z)=\omega^{2} d(z)$;\\
\rm{(iv)} Euler's formula
\begin{equation*}
\mathrm{e}^{\omega^{k}z}=c(z)+\omega^{k}s(z)+\frac{1}{\omega ^{k}}d(z),\quad k=1,2,3;
\end{equation*}
\rm{(v)} the main identity
\begin{equation*}
c^{3}(z)+s^{3}(z)+d^{3}(z)-3c(z)s(z)d(z)=1;
\end{equation*}
\rm{(vi)} the summation formulas
\begin{eqnarray*}
&&c(z_{1}+z_{2})=c(z_{1})c(z_{2})+s(z_{1})d(z_{2})+d(z_{1})s(z_{2}),\\
&&s(z_{1}+z_{2})=c(z_{1})s(z_{2})+s(z_{1})c(z_{2})+d(z_{1})d(z_{2}),\\
&&d(z_{1}+z_{2})=c(z_{1})d(z_{2})+s(z_{1})s(z_{2})+d(z_{1})c(z_{2});
\end{eqnarray*}
\rm{(vii)}
\begin{eqnarray*}
&&3c^{2}(z)=c(2z)+2c(-z),\\
&&3s^{2}(z)=d(2z)+2d(-z),\\
&&3d^{2}(z)=s(2z)+2s(-z);
\end{eqnarray*}
\rm{(viii)}
\begin{eqnarray*}
&&c^{2}(z)-s(z)d(z)=c(-z),\\
&&d^{2}(z)-s(z)c(z)=s(-z),\\
&&s^{2}(z)-d(z)c(z)=d(-z);
\end{eqnarray*}
\rm{(ix)} Taylor's formulas
\begin{eqnarray*}
&&c(z)=\sum_{n=0}^{\infty}\frac{z^{3n}}{3n!}=1+\frac{z^{3}}{3!}+\frac{z^{6}}{6!}+\cdots+\frac{z^{3n}}{(3n)!}+\cdots,\\
&&s(z)=\sum_{n=0}^{\infty}\frac{z^{3n+1}}{(3n+1)!}=z+\frac{z^{4}}{4!}+\frac{z^{7}}{7!}+\cdots+\frac{z^{3n+1}}{(3n+1)!}+\cdots,\\
&&d(z)=\sum_{n=0}^{\infty}\frac{z^{3n+2}}{(3n+2)!}=\frac{z^{2}}{2!}+\frac{z^{5}}{5!}+\frac{z^{8}}{8!}+\cdots+\frac{z^{3n+2}}{(3n+2)!}+\cdots.
\end{eqnarray*}
\end{lem}
\Proof
From Proposition \ref{omega}, relations in (i)-(ix) are obtained.
\qed

\begin{lem}
\label{csd2npi}
If $z={2n\pi}\mathrm{i}$, $n\in\mathbb{Z}$, then
\begin{eqnarray*}
&&c(2n\pi\mathrm{i})=\frac{1}{3}\left(1+2(-1)^{n}\cosh \sqrt{3}n\pi\right),\\
&&s(2n\pi\mathrm{i})=\frac{1}{3}\left(1-(-1)^{n}\cosh \sqrt{3}n\pi+\sqrt{3}\mathrm{i}(-1)^{n}\sinh \sqrt{3}n\pi\right),\\
&&d(2n\pi\mathrm{i})=\frac{1}{3}\left(1-(-1)^{n}\cosh \sqrt{3}n\pi-\sqrt{3}\mathrm{i}(-1)^{n}\sinh \sqrt{3}n\pi\right)
\end{eqnarray*}
and
\begin{equation*}
c(2n\pi\mathrm{i})+s(2n\pi\mathrm{i})+d(2n\pi\mathrm{i})=1.
\end{equation*}
Besides,
\begin{equation*}
c(-2n\pi\mathrm{i})=c(2n\pi\mathrm{i}),\quad d(-2n\pi\mathrm{i})=s(2n\pi\mathrm{i}),\quad s(-2n\pi\mathrm{i})=d(2n\pi\mathrm{i}).
\end{equation*}
\end{lem}
\Proof
According to the expressions (\ref{csd}) for $c(z)$, $s(z)$ and $d(z)$, a simple calculation will finish the proof.
\qed

Due to the expressions (\ref{csd}) of the functions $c(z)$, $s(z)$ and $d(z)$, the Fourier transforms for $v\in L_{\mathbb{C}}^{2}(0,1)$ are
\begin{eqnarray}
&&\widetilde{v}_{c}(\lambda):=\left<v(x),c(\mathrm{i}\overline{\lambda}x)\right>=\frac{1}{3}\sum\limits_{k=1}^{3}\widetilde{v}_{k}(\lambda),\notag\\
&&\widetilde{v}_{s}(\lambda):=\left<v(x),s(\mathrm{i}\overline{\lambda}x)\right>=\frac{1}{3}\sum\limits_{k=1}^{3}\frac{1}{\omega^{k}}\widetilde{v}_{k}(\lambda),\label{v}\\
&&\widetilde{v}_{d}(\lambda):=\left<v(x),d(\mathrm{i}\overline{\lambda}x)\right>=\frac{1}{3}\sum\limits_{k=1}^{3}\omega^{k}\widetilde{v}_{k}(\lambda),\notag
\end{eqnarray}
where
\begin{equation*}
\widetilde{v}_{k}(\lambda):=\int_{0}^{1}\mathrm{e}^{-\mathrm{i}\omega^{k}\lambda x}v(x)\mathrm{d} x,\quad k=1,2,3.
\end{equation*}
Let
\begin{equation*}
f^{*}(\lambda)=\overline{f(\overline{\lambda})}
\end{equation*}
denote the operation of involution, then
\begin{eqnarray}
&&\widetilde{v}_{c}^{*}(\lambda):=\left<c(\mathrm{i}{\lambda}x),v(x)\right>=\frac{1}{3}\sum\limits_{k=1}^{3}\widetilde{v}_{k}^{*}(\lambda),\notag\\
&&\widetilde{v}_{s}^{*}(\lambda):=\left<s(\mathrm{i}{\lambda}x),v(x)\right>=\frac{1}{3}\sum\limits_{k=1}^{3}\omega^{k}\widetilde{v}_{k}^{*}(\lambda),\label{vstar}\\
&&\widetilde{v}_{d}^{*}(\lambda):=\left<d(\mathrm{i}{\lambda}x),v(x)\right>=\frac{1}{3}\sum\limits_{k=1}^{3}\frac{1}{\omega^{k}}\widetilde{v}_{k}^{*}(\lambda),\notag
\end{eqnarray}
here,
\begin{equation*}
\widetilde{v}_{k}^{*}(\lambda):=\int_{0}^{1}\mathrm{e}^{\mathrm{i}\frac{1}{\omega^{k}}\lambda x}\overline{v}(x)\mathrm{d} x,\quad k=1,2,3.
\end{equation*}

\begin{rem}
For $k=1,2,3$, $\widetilde{v}_{k}(\lambda)$ and $\widetilde{v}_{k}^{*}(\lambda)$ are entire functions of exponential type, and hence $(\ref{v})$, $(\ref{vstar})$ are entire functions of exponential type.
\end{rem}
\begin{lem}
\label{wcsd}
For $v\in L_{\mathbb{C}}^{2}(0,1)$, the following equalities holds:
\begin{eqnarray*}
&&c(\mathrm{i}\lambda)\widetilde{v}_{d}(\lambda)+s(\mathrm{i}\lambda)\widetilde{v}_{s}(\lambda)+d(\mathrm{i}\lambda)\widetilde{v}_{c}(\lambda)
=\widetilde{w}_{d}(-\lambda),\\
&&d(\mathrm{i}\lambda)\widetilde{v}_{d}(\lambda)+c(\mathrm{i}\lambda)\widetilde{v}_{s}(\lambda)+s(\mathrm{i}\lambda)\widetilde{v}_{c}(\lambda)
=\widetilde{w}_{s}(-\lambda),\\
&&s(\mathrm{i}\lambda)\widetilde{v}_{d}(\lambda)+d(\mathrm{i}\lambda)\widetilde{v}_{s}(\lambda)+c(\mathrm{i}\lambda)\widetilde{v}_{c}(\lambda)
=\widetilde{w}_{c}(-\lambda),
\end{eqnarray*}
where
\begin{eqnarray*}
&&\widetilde{w}_{c}(\lambda)=\left<w(x),c(\mathrm{i}\overline{\lambda}x)\right>, \\
&&\widetilde{w}_{s}(\lambda)=\left<w(x),s(\mathrm{i}\overline{\lambda}x)\right>, \\
&&\widetilde{w}_{d}(\lambda)=\left<w(x),d(\mathrm{i}\overline{\lambda}x)\right>
\end{eqnarray*}
are the Fourier transforms for $w(x):=v(1-x)$.
\end{lem}
\Proof
It follows from the definitions of $\widetilde{v}_{c}(\lambda)$, $\widetilde{v}_{s}(\lambda)$, $\widetilde{v}_{d}(\lambda)$, and the relations (ii), (vi) in Lemma \ref{entire} that
\begin{eqnarray*}
c(\mathrm{i}\lambda)\widetilde{v}_{d}(\lambda)+s(\mathrm{i}\lambda)\widetilde{v}_{s}(\lambda)+d(\mathrm{i}\lambda)\widetilde{v}_{c}(\lambda)
\hh &=&\hh \int_{0}^{1}d(\mathrm{i}\lambda(1-t))v(t)\mathrm{d}t\\
\hh &=&\hh \int_{0}^{1}d(\mathrm{i}\lambda x)v(1-x)\mathrm{d}x\\
\hh &=&\hh \widetilde{w}_{d}(-\lambda).
\end{eqnarray*}
Similarly, we can proof the remaining equalities.
\qed

Set
\begin{equation}
m(\lambda):=\left<\int_{0}^{x}d(\mathrm{i}\lambda(x-t))v(t)\mathrm{d}t,v(x)\right>=\frac{1}{3}\sum\limits_{k=1}^{3}\omega^{k}\phi_{k}(\lambda),\label{mlambda}
\end{equation}
where
\begin{equation*}
\phi_{k}(\lambda):=\int_{0}^{1}\mathrm{d}x\int_{0}^{x}\mathrm{e}^{\mathrm{i}\omega^{k}\lambda (x-t)}\overline{v}(x){v}(t)\mathrm{d}t,\quad k=1,2,3
\end{equation*}
are the Fourier transforms of the convolution.
\begin{lem}
\label{mmstar}
The function $m(\lambda)$ is an entire function of exponential type and satisfies
\begin{equation*}
m(\lambda)+m^{*}(\lambda)=\widetilde{v}_{d}(\lambda)\widetilde{v}_{c}^{*}(\lambda)+\widetilde{v}_{s}(\lambda)\widetilde{v}_{s}^{*}(\lambda)+\widetilde{v}_{c}(\lambda)\widetilde{v}_{d}^{*}(\lambda).
\end{equation*}
\end{lem}
\Proof
See \cite[Lemma 2.1]{VAZ2021}.
\qed

\section{The operator $L_{0}$}
In this section, we will investigate the spectrum, eigenfunctions and resolvent of the operator $L_0$, which are fundamental for discussing the relevant properties of the operator $L_{\alpha}$. The domain of the operator $L_0$ is
\begin{equation*}
\mathcal{D}(L_{0})=\left\{y\in W_3^2 (0,1)\vert y(1)=y(0), y^{\prime}(1)=y^{\prime}(0), y^{\prime\prime}(1)=y^{\prime\prime}(0)\right\}.
\end{equation*}
By (\ref{csd}), we find $c(\mathrm{i}\lambda x)$, $\frac{s(\mathrm{i}\lambda x)}{\mathrm{i}\lambda}$ and $\frac{d(\mathrm{i}\lambda x)}{(\mathrm{i}\lambda)^{2}}$ are the solutions of
\begin{equation}
\mathrm{i}y^{\prime\prime\prime}(x)=zy(x),\quad z=\lambda^{3},\quad \lambda \in \mathbb{C} \label{L0}
\end{equation}
satisfying the initial conditions
\begin{equation*}
\left(
\begin{array}{ccc}
c(0) & \frac{s(0)}{\mathrm{i}\lambda} & \frac{d(0)}{(\mathrm{i}\lambda)^2} \\
\mathrm{i}\lambda c^{\prime }(0) & s^{\prime }(0) &\frac{d^{\prime }(0)}{\mathrm{i}\lambda} \\
 (\mathrm{i}\lambda)^2 c^{\prime \prime}(0) &\mathrm{i}\lambda s^{\prime \prime}(0) & d^{\prime \prime}(0)%
\end{array}%
\right) =I_{3}.
\end{equation*}
Then each unique solution of (\ref{L0}) with initial conditions
\begin{equation}
\left(y(0),y^{\prime}(0),y^{\prime\prime}(0)\right)=\left(a_{0},a_{1},a_{2}\right)\in \mathbb{C}^3 \label{initial}
\end{equation}
is
\begin{equation}
y_{0}(x,\lambda)=a_{0}c(\mathrm{i}\lambda x)+a_{1}\frac{s(\mathrm{i}\lambda x)}{\mathrm{i}\lambda}+a_{2}\frac{d(\mathrm{i}\lambda x)}{(\mathrm{i}\lambda)^{2}}. \label{solutionl0}
\end{equation}
Thus, using the method of variation of constants, one finds the unique solution to the third order inhomogeneous differential equation
\begin{equation}
\mathrm{i}y^{\prime\prime\prime}(x)=\lambda^{3}y(x)+f(x), \quad f\in L_{\mathbb{C}}^2(0,1)\label{cauchy}
\end{equation}
satisfying the initial conditions (\ref{initial}) is
\begin{equation}
y(x,\lambda)=a_{0}c(\mathrm{i}\lambda x)+a_{1}\frac{s(\mathrm{i}\lambda x)}{\mathrm{i}\lambda}+a_{2}\frac{d(\mathrm{i}\lambda x)}{(\mathrm{i}\lambda)^{2}}+\mathrm{i}\int_{0}^{x}\frac{d(\mathrm{i}\lambda(x-t))}{\lambda^{2}}f(t)\mathrm{d}t.\label{cauchys}
\end{equation}
\subsection{The characteristic function, spectrum and eigenfunctions of the operator $L_{0}$}

In this subsection, we give the definition and properties for the characteristic function, the spectrum and the eigenfunctions of $L_{0}$.
\begin{lem}
\label{l0s}
Let
\begin{equation*}
M(0,\lambda):=\left(
\begin{array}{ccc}
c(\mathrm{i}\lambda )-1 & {\frac{s(\mathrm{i}\lambda l)}{\mathrm{i}\lambda}} &  {\frac{d(\mathrm{i}\lambda )}{(\mathrm{i}\lambda)^{2}}} \\
(\mathrm{i}\lambda)c^{\prime }(\mathrm{i}\lambda ) & s^{\prime }(\mathrm{i}\lambda )-1 & {\frac{d^{\prime }(\mathrm{i}\lambda )}{\mathrm{i}\lambda}} \\
 (\mathrm{i}\lambda)^{2}c^{\prime \prime}(\mathrm{i}\lambda ) &(\mathrm{i}\lambda)s^{\prime \prime}(\mathrm{i}\lambda ) & d^{\prime \prime }(\mathrm{i}\lambda )-1%
\end{array}%
\right).
\end{equation*}
The spectrum of the operator $L_{0}$ is denoted by
\begin{equation*}
\sigma(L_{0}):=\{z_{n}:=\lambda_{n}^{3}(0)\vert\lambda_n={2n\pi},n\in \mathbb{Z}\},
\end{equation*}
where $\lambda_{n}(0)$ is the root of the characteristic function
\begin{equation}
\Delta (0,\lambda):=\det M(0,\lambda)=3[c(\mathrm{i}\lambda )-c(-\mathrm{i}\lambda )]=-8\mathrm{i}\sin\frac{\lambda }{2}\sin\frac{\omega\lambda }{2}\sin\frac{\omega^{2}\lambda }{2}=0. \label{Delta0}
\end{equation}
Besides, the characteristic function $\Delta (0,\lambda)$ satisfies
\begin{equation}
\Delta (0,\lambda)=\Delta (0,\omega\lambda),\quad \Delta ^{*}(0,\lambda)=-\Delta (0,\lambda).\label{d0pp}
\end{equation}
\end{lem}
\Proof
Since the unique solution of (\ref{L0})-(\ref{initial}) is given in (\ref{solutionl0}), then it follows from the periodic boundary conditions (\ref{bc}) that the system of linear equations for $a_0$, $a_1$, $a_2$,
\begin{equation}
M(0,\lambda)\left(
\begin{array}{c}
a_{0} \\
a_{1}\\
a_{2}
\end{array}
\right)=\left(\begin{array}{c}
0 \\
0\\
0
\end{array}%
\right).\label{ma0a1a2}
\end{equation}
Obviously, $z=\lambda^{3}$ is an eigenvalue of the operator $L_{0}$ if and only if the system (\ref{ma0a1a2}) has non-trivial solution, i.e., $\lambda$ is the root of $\Delta (0,\lambda)=\det M(0,\lambda)=0$. According to Proposition \ref{omega} and (v), (viii) in Lemma \ref{entire}, we can simplify the characteristic function
\begin{equation}
\Delta (0,\lambda)=3[c(\mathrm{i}\lambda )-c(-\mathrm{i}\lambda )]\label{3cc}
\end{equation}
and obtain the equality
\begin{equation*}
c(\mathrm{i}\lambda x)=\frac{1}{3}\left(4\cos\frac{\lambda x}{2}\cos\frac{\omega\lambda x}{2}\cos\frac{\omega^2\lambda x}{2}-1-4\mathrm{i}\sin\frac{\lambda x}{2}\sin\frac{\omega\lambda x}{2}\sin\frac{\omega^2\lambda x}{2}\right).
\end{equation*}
Hence, the identity (\ref{Delta0}) holds and the roots of $\Delta (0,\lambda)=0$ are $\lambda_{n}(0)={2n\pi}$, $n\in \mathbb{Z}$. Furthermore, due to (ii), (iii) in Lemma \ref{entire}, the relations in (\ref{d0pp}) are also proved.\qed

Now we deduce the multiplicity of the each eigenvalue and the expressions of the eigenfunctions.
\begin{lem}
\label{mef}
\rm{(i)} The eigenvalues of the operator $L_{0}$ are $z_n=(2n\pi)^3$ and the multiplicity of each eigenvalue $z_{n}$ equals to one. \\
\rm{(ii)} The eigenfunctions of the operator corresponding to $z_n$ are
\begin{equation*}
u_n(0,x)=\mathrm{e}^{2n\pi\mathrm{i}x}, \quad n\in\mathbb{Z}.
\end{equation*}
\end{lem}
\Proof
\rm{(i)} Due to (i) in Lemma \ref{entire}, we see
\begin{eqnarray*}
M(0,\lambda_{n}(0))\hh &=&\hh \left(
\begin{array}{ccc}
c(\mathrm{i}\lambda_{n}(0))-1 &  {\frac{s(\mathrm{i}\lambda_{n}(0))}{\mathrm{i}\lambda_{n}(0)} }& {\frac{d(\mathrm{i}\lambda_{n}(0))}{(\mathrm{i}\lambda_{n}(0))^{2}} }\\
\mathrm{i}\lambda_{n}(0)c^{\prime }(\mathrm{i}\lambda_{n}(0)) & s^{\prime }(\mathrm{i}\lambda_{n}(0))-1 & {\frac{d^{\prime }(\mathrm{i}\lambda_{n}(0))}{\mathrm{i}\lambda_{n}(0)} }\\
 (\mathrm{i}\lambda_{n}(0))^{2}c^{\prime \prime}(\mathrm{i}\lambda_{n}(0)) &\mathrm{i}\lambda_{n}(0)s^{\prime \prime}(\mathrm{i}\lambda_{n}(0)) & d^{\prime \prime }(\mathrm{i}\lambda_{n}(0))-1%
\end{array}%
\right)\\
\hh &=&\hh \left(
\begin{array}{ccc}
c(\mathrm{i}\lambda_{n}(0))-1 &  {\frac{s(\mathrm{i}\lambda_{n}(0))}{\mathrm{i}\lambda_{n}(0)} }& {\frac{d(\mathrm{i}\lambda_{n}(0))}{(\mathrm{i}\lambda_{n}(0))^{2}}} \\
\mathrm{i}\lambda_{n}(0)d(\mathrm{i}\lambda_{n}(0)) & c(\mathrm{i}\lambda_{n}(0))-1 & {\frac{s(\mathrm{i}\lambda_{n}(0))}{\mathrm{i}\lambda_{n}(0)}} \\
 (\mathrm{i}\lambda_{n}(0))^{2}s(\mathrm{i}\lambda_{n}(0)) &\mathrm{i}\lambda_{n}(0)d(\mathrm{i}\lambda_{n}(0)) & c(\mathrm{i}\lambda_{n}(0))-1%
\end{array}%
\right).
\end{eqnarray*}

For $z_{0}=0$, it follows from
\begin{equation*}
c(0)=1,\quad\lim_{\lambda\rightarrow 0}\frac{s(\mathrm{i}\lambda )}{\mathrm{i}\lambda}=1,\quad \lim_{\lambda\rightarrow 0}\frac{d(\mathrm{i}\lambda )}{(\mathrm{i}\lambda)^{2}}=\frac{1}{2}
\end{equation*}
that
\begin{equation*}
M(0,0)= \left(
\begin{array}{ccc}
0 & 1 & \frac{1}{2} \\
0 & 0 &1 \\
0 &0 & 0%
\end{array}%
\right)\sim\left(
\begin{array}{ccc}
0 & 1 & 0 \\
0 & 0 &1 \\
0 &0 & 0%
\end{array}%
\right).
\end{equation*}
Here, the notation '$\sim$' means the equivalence relation in matrix theory. Hence, $\text{Rank}M(0,0)=2$ and $z_{0}=0$ is a simple eigenvalue of the operator $L_{0}$.

For any $z_{n}=({2n\pi})^3$, $n\in\mathbb{Z}/\{0\}$, according to Lemma \ref{csd2npi}, we see $c(2n\pi\mathrm{i})-1\neq 0$ holds for all $n\neq 0$. Then, from (viii) in Lemma \ref{entire} and Lemma \ref{csd2npi}, a straightforward calculation yields that
\begin{eqnarray}
M(0,2{n}\pi)\hh &\sim &\hh \left(
\begin{array}{ccc}
c(2n\pi\mathrm{i})-1 &  {\frac{s(2n\pi\mathrm{i})}{2n\pi\mathrm{i}}} &  {\frac{d(2n\pi\mathrm{i})}{(2n\pi\mathrm{i})^{2}} }\\
2n\pi\mathrm{i}d(2n\pi\mathrm{i}) & c(2n\pi\mathrm{i})-1 & {\frac{s(2n\pi\mathrm{i})}{2n\pi\mathrm{i}}} \\
0 &0 & 0%
\end{array}%
\right)\notag\\
\hh &\sim &\hh \left(
\begin{array}{ccc}
 1&  {\frac{s(2n\pi\mathrm{i})}{2n\pi\mathrm{i}[c(2n\pi\mathrm{i})-1]}} & {\frac{d(2n\pi\mathrm{i})}{(2n\pi\mathrm{i})^{2}[c(2n\pi\mathrm{i})-1]}} \\
0 &  {\frac{c(2n\pi\mathrm{i})-1}{2n\pi\mathrm{i}}-\frac{s(2n\pi\mathrm{i})d(2n\pi\mathrm{i})}{2n\pi\mathrm{i}[c(2n\pi\mathrm{i})-1]}} & {\frac{s(2n\pi\mathrm{i})}{(2n\pi\mathrm{i})^{2}}-\frac{d^{2}(2n\pi\mathrm{i})}{(2n\pi\mathrm{i})^{2}[c(2n\pi\mathrm{i})-1]} }\\
0 &0 & 0%
\end{array}%
\right)\notag\\
\hh &=&\hh \left(
\begin{array}{ccc}
 1&  {\frac{s(2n\pi\mathrm{i})}{2n\pi\mathrm{i}[c(2n\pi\mathrm{i})-1]}} & {\frac{d(2n\pi\mathrm{i})}{(2n\pi\mathrm{i})^{2}[c(2n\pi\mathrm{i})-1]}} \\
0& {-\frac{1}{2n\pi\mathrm{i}}}& {\frac{1}{(2n\pi\mathrm{i})^{2}}}\\
0 &0 & 0%
\end{array}%
\right)\notag\\
\hh &\sim&\hh \left(
\begin{array}{ccc}
 1& 0 &  {-\frac{1}{(2n\pi\mathrm{i})^{2}}} \\
0 & 1& {-\frac{1}{2n\pi\mathrm{i}}} \\
0 &0 & 0%
\end{array}%
\right).\label{eqv}
\end{eqnarray}
Hence, $\text{Rank}M(0,2n\pi)=2$. This implies that the multiplicity of $z_{n}=({2n\pi})^3$, $n\in\mathbb{Z}/\{0\}$ is one.

\rm{(ii)} For $z_{0}=0$, we find that the non-trivial solution of equation $\mathrm{i}y^{\prime\prime\prime}(x)=0$ satisfying the periodic boundary conditions (\ref{bc}) is $y_{0}(x,0)=C$, here $C$ denotes a nonzero constant. So the normalized eigenfunction with respect to $z_{n}=0$ is $u_0(0,x)=\mathrm{e}^{0}$.

For any $z_{n}=({2n\pi})^3$, $n\in\mathbb{Z}/\{0\}$, it follows from (\ref{eqv}) that the non-trivial solution of equation $\mathrm{i}y^{\prime\prime\prime}(x)=(2n\pi)^{3}y(x)$ determined by the periodic boundary conditions (\ref{bc}) is $y_{0}(x,2n\pi\mathrm{i})=c(2n\pi\mathrm{i}x)+s(2n\pi\mathrm{i}x)+d(2n\pi\mathrm{i}x)$. Hence the normalized eigenfunction with respect to $z_{n}=({2n\pi})^3$, $n\in\mathbb{Z}/\{0\}$ is $u_n(0,x)=\mathrm{e}^{2n\pi\mathrm{i}x}$.
\qed
\begin{rem}
The eigenfunctions $u_{n}:=u_n(0,x)$ of the operator $L_{0}$ satisfy
\begin{equation*}
\left<u_{n},u_{m}\right>=\delta_{n,m}:=\left\{\begin{array}{ll}
1,&m=n,\\
0,&m\neq n,
\end{array}
\right.\quad m,n\in \mathbb{Z}
\end{equation*}
and form an orthonormal basis in the Hilbert space $L_{\mathbb{C}}^{2}(0,1)$.
\end{rem}
\subsection{The resolvent of the operator $L_{0}$}
In this subsection, we give the expressions for the resolvent $R_{L_{0}}$ of the operator $L_{0}$, which play an important role in the inverse problems of the operator $L_{\alpha}$.
\begin{lem}
\rm{\cite{kato}} The resolvent $R_{L_{0}}(z)=(L_{0}-zI)^{-1}$ of the operator $L_{0}$ is
\begin{equation}
R_{L_{0}}(z)f=\sum_{n\in \mathbb{Z}}\frac{f_{n}}{z_{n}-z}u_{n},\label{rl0}
\end{equation}
where $f_{n}:=\left<f,u_{n}\right>$ are Fourier coefficients of $f\in L_{\mathbb{C}}^{2}(0,1)$ in the basis $\{u_{n}\}$.
\end{lem}

\begin{lem}
\label{rel0}
For $z=\lambda^3$, $f\in L_{\mathbb{C}}^{2}(0,1)$, the resolvent $R_{L_{0}}$ of the operator $L_{0}$ can be expressed by
\begin{eqnarray*}
(R_{L_{0}}(\lambda^{3})f)(x)
\hh &=&\hh \frac{-\mathrm{i}}{\lambda^{2}\Delta(0,\lambda)}\int_{0}^{x}\left\{d(\mathrm{i}\lambda(x-t))[1-3c(\mathrm{i}\lambda)]+d(-\mathrm{i}\lambda t)c(\mathrm{i}\lambda(x-1))\right.\\
\hh &&\hh \left.+s(-\mathrm{i}\lambda t)s(\mathrm{i}\lambda(x-1))+c(-\mathrm{i}\lambda t)d(\mathrm{i}\lambda(x-1))+d(\mathrm{i}\lambda x)c(\mathrm{i}\lambda(1-t))\right.\\
\hh &&\hh \left.+s(\mathrm{i}\lambda x)s(\mathrm{i}\lambda(1-t))+c(\mathrm{i}\lambda x)d(\mathrm{i}\lambda(1-t))\right\}f(t)\mathrm{d}t\\
\hh &&\hh +\frac{-\mathrm{i}}{\lambda^{2}\Delta(0,\lambda)}\int_{x}^{1}\left\{d(\mathrm{i}\lambda(x-t))[1-3c(-\mathrm{i}\lambda)]+d(-\mathrm{i}\lambda t)c(\mathrm{i}\lambda(x-1))\right.\\
\hh &&\hh \left.+s(-\mathrm{i}\lambda t)s(\mathrm{i}\lambda(x-1))+c(-\mathrm{i}\lambda t)d(\mathrm{i}\lambda(x-1))+d(\mathrm{i}\lambda x)c(\mathrm{i}\lambda(1-t))\right.\\
\hh &&\hh \left.+s(\mathrm{i}\lambda x)s(\mathrm{i}\lambda(1-t))+c(\mathrm{i}\lambda x)d(\mathrm{i}\lambda(1-t))\right\}f(t)\mathrm{d}t.
\end{eqnarray*}
\end{lem}
\Proof
Let $R_{L_{0}}(\lambda^{3})f=(L_{0}-\lambda^3I)^{-1}f=y$, then $L_{0}y=\lambda^{3}y+f$ and $(\ref{cauchys})$ is the solution of it satisfying the initial conditions (\ref{initial}). Note that the periodic boundary conditions (\ref{bc}) for $y(x,\lambda)$ in (\ref{cauchys}) lead to the system of linear equations
\begin{equation*}
M(0,\lambda)\left(
\begin{array}{c}
a_{0} \\
a_{1}\\
a_{2}
\end{array}
\right)=\left(\begin{array}{c}
-\mathrm{i} \int_{0}^{1}\frac{d(\mathrm{i}\lambda(1-t))}{\lambda^{2}}f(t)\mathrm{d}t \\
 \int_{0}^{1}\frac{s(\mathrm{i}\lambda(1-t))}{\lambda}f(t)\mathrm{d}t \\
\mathrm{i} \int_{0}^{1}c(\mathrm{i}\lambda(1-t))f(t)\mathrm{d}t
\end{array}%
\right)
\end{equation*}
relative to $a_{0}$, $a_{1}$, $a_{2}$, where $M(0,\lambda)$ is given in Lemma \ref{l0s}.
%Since
%\begin{eqnarray*}
%\hh \hh &&\hh \hh M^{-1}(0,\lambda)\\
%\hh \hh &=&\hh \hh \frac{1}{\Delta(0,\lambda)}\left(\hh
%\begin{array}{ccc}
%[c(\mathrm{i}\lambda)-1]^{2}-s(\mathrm{i}\lambda)d(\mathrm{i}\lambda) & %\frac{1}{\mathrm{i}\lambda}[d^{2}(\mathrm{i}\lambda)-s(\mathrm{i}\lambda)(c(\mathrm{i}\lambda)-1)] & \frac{1}{(\mathrm{i}\lambda)^{2}}[s^{2}(\mathrm{i}\lambda)-d(\mathrm{i}\lambda)(c(\mathrm{i}\lambda)-1)] \\
%\mathrm{i}\lambda[s^{2}(\mathrm{i}\lambda)-d(\mathrm{i}\lambda)(c(\mathrm{i}\lambda)-1)]  & [c(\mathrm{i}\lambda)-1]^{2}-s(\mathrm{i}\lambda)d(\mathrm{i}\lambda)&\frac{1}{\mathrm{i}\lambda}[d^{2}(\mathrm{i}\lambda)-s(\mathrm{i}\lambda)(c(\mathrm{i}\lambda)-1)] \\
%{(\mathrm{i}\lambda)^{2}}[d^{2}(\mathrm{i}\lambda)-s(\mathrm{i}\lambda)(c(\mathrm{i}\lambda)-1)] &\mathrm{i}\lambda[s^{2}(\mathrm{i}\lambda)-d(\mathrm{i}\lambda)(c(\mathrm{i}\lambda)-1)]& [c(\mathrm{i}\lambda)-1]^{2}-s(\mathrm{i}\lambda)d(\mathrm{i}\lambda)%
%\end{array}%
%\hh \right)
%\end{eqnarray*}
Then, due to (v), (vi) and (viii) in Lemma \ref{entire}, a straightforward calculation gives
\begin{eqnarray*}
a_{0}\hh &=&\hh \frac{-\mathrm{i}}{\lambda^{2}{\Delta(0,\lambda)}}\int_{0}^{1}\left\{\left[(c(\mathrm{i}\lambda)-1)^{2}-s(\mathrm{i}\lambda)d(\mathrm{i}\lambda)\right]d(\mathrm{i}\lambda(1-t))-\left[s(\mathrm{i}\lambda)(c(\mathrm{i}\lambda)-1)\right.\right.\\
\hh &&\hh \left.\left.-d^{2}(\mathrm{i}\lambda)\right]s(\mathrm{i}\lambda(1-t))+[s^{2}(\mathrm{i}\lambda)-d(\mathrm{i}\lambda)(c(\mathrm{i}\lambda)-1)]c(\mathrm{i}\lambda(1-t))\right\}f(t)\mathrm{d}t\\
\hh &=&\hh \frac{-\mathrm{i}}{\lambda^{2}{\Delta(0,\lambda)}}\int_{0}^{1}\left\{d(-\mathrm{i}\lambda t)[1+c(\mathrm{i}\lambda)-2c(-\mathrm{i}\lambda)]+s(-\mathrm{i}\lambda t)[s(\mathrm{i}\lambda)+s(-\mathrm{i}\lambda)]\right.\\
\hh &&\hh \left.+c(-\mathrm{i}\lambda t)[d(\mathrm{i}\lambda)+d(-\mathrm{i}\lambda)]\right\}f(t)\mathrm{d}t,
\end{eqnarray*}
\begin{eqnarray*}
a_{1}\hh &=&\hh \frac{1}{\lambda{\Delta(0,\lambda)}}\int_{0}^{1}\left\{[s^{2}(\mathrm{i}\lambda)-d(\mathrm{i}\lambda)(c(\mathrm{i}\lambda)-1)]d(\mathrm{i}\lambda(1-t))+\left[(c(\mathrm{i}\lambda)-1)^{2}\right.\right.\\
\hh &&\hh \left.\left.-s(\mathrm{i}\lambda)d(\mathrm{i}\lambda)\right]s(\mathrm{i}\lambda(1-t))+[d^{2}(\mathrm{i}\lambda)-s(\mathrm{i}\lambda)(c(\mathrm{i}\lambda)-1)]c(\mathrm{i}\lambda(1-t))\right\}f(t)\mathrm{d}t\\
\hh &=&\hh \frac{1}{\lambda{\Delta(0,\lambda)}}\int_{0}^{1}\left\{d(-\mathrm{i}\lambda t)[d(\mathrm{i}\lambda)+d(-\mathrm{i}\lambda)]+s(-\mathrm{i}\lambda t)[1+c(\mathrm{i}\lambda)-2c(-\mathrm{i}\lambda)]\right.\\
\hh &&\hh \left.+c(-\mathrm{i}\lambda t)[s(\mathrm{i}\lambda)+s(-\mathrm{i}\lambda)]\right\}f(t)\mathrm{d}t,\\
a_{2}\hh &=&\hh \frac{1}{{\Delta(0,\lambda)}}\int_{0}^{1}\left\{[d^{2}(\mathrm{i}\lambda)-s(\mathrm{i}\lambda)(c(\mathrm{i}\lambda)-1)]d(\mathrm{i}\lambda(1-t))-\left[d(\mathrm{i}\lambda)(c(\mathrm{i}\lambda)-1)\right.\right.\\
\hh &&\hh \left.\left.-s^{2}(\mathrm{i}\lambda)\right]s(\mathrm{i}\lambda(1-t))+[(c(\mathrm{i}\lambda)-1)^{2}-s(\mathrm{i}\lambda)d(\mathrm{i}\lambda)]c(\mathrm{i}\lambda(1-t))\right\}f(t)\mathrm{d}t\\
\hh &=&\hh \frac{1}{{\Delta(0,\lambda)}}\int_{0}^{1}\left\{d(-\mathrm{i}\lambda t)[s(\mathrm{i}\lambda)+s(-\mathrm{i}\lambda)]+s(-\mathrm{i}\lambda t)[d(\mathrm{i}\lambda)+d(-\mathrm{i}\lambda)]\right.\\
\hh &&\hh \left.+c(-\mathrm{i}\lambda t)[1+c(\mathrm{i}\lambda)-2c(-\mathrm{i}\lambda)]\right\}f(t)\mathrm{d}t.
\end{eqnarray*}
Therefore, it follows from (vi) in Lemma \ref{entire} that
\begin{eqnarray*}
y(x,\lambda)\hh &=&\hh a_{0}c(\mathrm{i}\lambda x)+a_{1}\frac{s(\mathrm{i}\lambda x)}{\mathrm{i}\lambda}+a_{2}\frac{d(\mathrm{i}\lambda x)}{(\mathrm{i}\lambda)^{2}}+\mathrm{i}\int_{0}^{x}\frac{d(\mathrm{i}\lambda(x-t))}{\lambda^{2}}f(t)\mathrm{d}t\\
\hh &=&\hh \frac{-\mathrm{i}}{\lambda^{2}\Delta (0,\lambda)}\int_{0}^{1}\left\{d(\mathrm{i}\lambda({x}- t))[1-3c(-\mathrm{i}\lambda)]+d(-\mathrm{i}\lambda t)c(\mathrm{i}\lambda(x-1))\right.\\
\hh &&\hh \left.+s(-\mathrm{i}\lambda t)s(\mathrm{i}\lambda(x-1))+c(-\mathrm{i}\lambda t)d(\mathrm{i}\lambda(x-1))+d(\mathrm{i}\lambda{x})c(\mathrm{i}\lambda(1-t))\right.\\
\hh &&\hh \left.+s(\mathrm{i}\lambda{x})s(\mathrm{i}\lambda(1-t))+c(\mathrm{i}\lambda{x})d(\mathrm{i}\lambda(1-t))\right\}f(t)\mathrm{d}t +\frac{\mathrm{i}}{\lambda^{2}}\int_{0}^{x}d(\mathrm{i}\lambda(1-t))f(t)\mathrm{d}t.
\end{eqnarray*}
Finally, on account of the equality (\ref{3cc}), we can prove this lemma.
\qed

\section{The operator $L_{\alpha}$}
In this section, we discuss the spectrum, eigenfunctions and resolvent for the operator $L_{\alpha}=L_{\alpha}(v)$
\begin{eqnarray*}
&&(L_{\alpha}y)(x)=\mathrm{i}y^{\prime\prime\prime}(x)+\alpha\int_{0}^{1}y(t)\overline{v}(t)\mathrm{d}tv(x),\quad\alpha\in\mathbb{R},\quad v\in L_{\mathbb{c}}^{2}(0,1),\\
&&\mathcal{D}(L_{\alpha})=\left\{y\in W_3^2 (0,1)\vert y(1)=y(0), y^{\prime}(1)=y^{\prime}(0), y^{\prime\prime}(1)=y^{\prime\prime}(0)\right\},
\end{eqnarray*}
which can be considered as a one-dimensional perturbation of the operator $L_{0}$.
%Without loss in generality, assume that $\Vert v\Vert_{L^{2}}=1$. Domains of the operators $L_{\alpha}$ and $L_{0}$ coincide, $\mathcal{D}(L_{\alpha})=\mathcal{D}(L_{0})$.
\subsection{The resolvent of the operator $L_{\alpha}$}
In this subsection, we show the expression for the resolvent of the operator $L_{\alpha}$ at first, which is essential for investigating the spectrum of the operator $L_{\alpha}$.
\begin{lem}
\label{rlaa}
The resolvent $R_{L_{\alpha}}(z)=(L_{\alpha}-zI)^{-1}$ of the operator $L_{\alpha}$ is expressed by the resolvent $R_{L_{0}}(z)=(L_{0}-zI)^{-1}$ of the operator $L_{0}$ and satisfies
\begin{eqnarray}
R_{L_{\alpha}}(z)f\hh &=&\hh R_{L_{0}}(z)f-\alpha\frac{\left<R_{L_{0}}(z)f,v\right>}{1+\alpha\left<R_{L_{0}}(z)v,v\right>}\cdot R_{L_{0}}(z)v\notag\\
\hh &=&\hh \sum\limits_{n\in\mathbb{Z}}\frac{f_{n}}{z_{n}-z}u_{n}-\alpha\frac{ \sum\limits_{k\in\mathbb{Z}}\frac{f_{k}\overline{v}_{k}}{z_{k}-z}}{1+\alpha \sum\limits_{k\in\mathbb{Z}}\frac{\vert v_{k}\vert^{2}}{z_{k}-z}}\cdot \sum\limits_{n\in\mathbb{Z}}\frac{v_{n}}{z_{n}-z}u_{n},\label{rla}
\end{eqnarray}
where $f_{n}=\left<f,u_{n}\right>$, $v_{n}=\left<v,u_{n}\right>$ are Fourier coefficients of $f,v\in L_{\mathbb{C}}^2(0,1)$  in the basis $\{u_{n}\}$ correspondingly.
\end{lem}
\Proof
Note that $\mathcal{D}(L_{\alpha})=\mathcal{D}(L_{0})$ and $L_{\alpha}y=L_{0}y+\alpha\left<y,v\right>v$. Denote $R_{L_{\alpha}}(z)f:=(L_{\alpha}-zI)^{-1}f=g$. Then it yields
\begin{equation*}
f=(L_{\alpha}-zI)g=(L_{0}-zI)g+\alpha\left<g,v\right>v,
\end{equation*}
and
\begin{equation}
g+\alpha\left<g,v\right>R_{L_{0}}(z)v=R_{L_{0}}(z)f,\quad R_{L_{0}}(z)=(L_{0}-zI)^{-1}.\label{gvv}
\end{equation}
Multiplying (\ref{gvv}) by $\overline{v}(x)$ and integrating it from $0$ to $1$, we find
\begin{equation*}
\left<g,v\right>+\alpha\left<g,v\right>\left<R_{L_{0}}(z)v,v\right>=\left<R_{L_{0}}(z)f,v\right>
\end{equation*}
and
\begin{equation*}
\left<g,v\right>=\frac{\left<R_{L_{0}}(z)f,v\right>}{1+\alpha\left<R_{L_{0}}(z)v,v\right>}.
\end{equation*}
Therefore,
\begin{equation*}
g=R_{L_{\alpha}}(z)f=R_{L_{0}}(z)f-\alpha\frac{\left<R_{L_{0}}(z)f,v\right>}{1+\alpha\left<R_{L_{0}}(z)v,v\right>}R_{L_{0}}(z)v.
\end{equation*}
If we plug (\ref{rl0}) back into the above equation, (\ref{rla}) is proved.%Substituting here (\ref{rl0}), we can prove (\ref{rla}).
\qed
\subsection{The characteristic function, spectrum and eigenfunctions of the operator $L_{\alpha}$}
According to (\ref{cauchy})-(\ref{cauchys}), the unique solution of $L_{\alpha}y=zy=\lambda^3 y$ satisfying the initial conditions (\ref{initial}) is
\begin{equation}
y(x,\lambda)=a_{0}c(\mathrm{i}\lambda x)+a_{1}\frac{s(\mathrm{i}\lambda x)}{\mathrm{i}\lambda}+a_{2}\frac{d(\mathrm{i}\lambda x)}{(\mathrm{i}\lambda)^{2}}-\mathrm{i}\alpha\left<y,v\right>\int_{0}^{x}\frac{d(\mathrm{i}\lambda(x-t))}{\lambda^{2}}v(t)\mathrm{d}t.\label{cauchysv}
\end{equation}
In this subsection, we show the expression for characteristic function of the operator $L_{\alpha}$ via (\ref{cauchysv}), and then discuss the spectrum of the operator $L_{\alpha}$ by using the expression of the resolvent $R_{L_{\alpha}}$. Finally, we explore the eigenfunctions of the operator $L_{\alpha}$.
\begin{lem}
Denote
\begin{equation*}
M(\alpha, \lambda):=
\left(\begin{array}{cccc}
\widetilde{v}_{c}^{*}(\lambda)& \frac{\widetilde{v}_{s}^{*}(\lambda)}{\mathrm{i}\lambda}
& \frac{\widetilde{v}_{d}^{*}(\lambda)}{(\mathrm{i}\lambda)^{2}}& {-(1+\frac{\mathrm{i}\alpha}{\lambda^2}m(\lambda))}\\
c(\mathrm{i}\lambda)-1& \frac{s(\mathrm{i}\lambda)}{\mathrm{i}\lambda}& \frac{d(\mathrm{i}\lambda)}{(\mathrm{i}\lambda)^{2}}
&-\mathrm{i}\alpha \int_{0}^{1}\frac{d(\mathrm{i}\lambda(1-t))}{\lambda^{2}}v(t)\mathrm{d}t\\
\mathrm{i}\lambda d(\mathrm{i}\lambda)&c(\mathrm{i}\lambda)-1& \frac{s(\mathrm{i}\lambda)}{\mathrm{i}\lambda}
&\alpha \int_{0}^{1}\frac{s(\mathrm{i}\lambda(1-t))}{\lambda}v(t)\mathrm{d}t\\
(\mathrm{i}\lambda)^{2}s(\mathrm{i}\lambda)&(\mathrm{i}\lambda)d(\mathrm{i}\lambda)&c(\mathrm{i}\lambda)-1
&\mathrm{i}\alpha \int_{0}^{1}c(\mathrm{i}\lambda(1-t))v(t)\mathrm{d}t
\end{array}
\right).
\end{equation*}
The characteristic function of the operator $L_{\alpha}$ is
\begin{equation}
\Delta(\alpha,\lambda):=\det M(\alpha, \lambda)=\Delta(0,\lambda)+\frac{\mathrm{i}\alpha}{\lambda^2}[F(\lambda)+F^{*}(\lambda)] \label{deltaa}
\end{equation}
where
\begin{equation}
F(\lambda):=m(\lambda)[3c(\mathrm{i}\lambda)-1]-\widetilde{v}_{c}(\lambda)\widetilde{w}_{d}^{*}(-\lambda)
-\widetilde{v}_{s}(\lambda)\widetilde{w}_{s}^{*}(-\lambda)-\widetilde{v}_{d}(\lambda)\widetilde{w}_{c}^{*}(-\lambda).\label{flambda}
\end{equation}
Besides, the characteristic function of the operator $L_{\alpha}$ satisfies
\begin{equation}
\Delta (\alpha,\lambda)=\Delta (\alpha,\omega\lambda),\quad \Delta ^{*}(\alpha,\lambda)=-\Delta (\alpha,\lambda).\label{dapp}
\end{equation}
\end{lem}
\Proof
Multiplying equality (\ref{cauchysv}) by $\overline{v}(x)$ and integrating it on $[0,1]$, we find
\begin{equation*}
a_{0}\widetilde{v}_{c}^{*}(\lambda)+a_{1}\frac{\widetilde{v}_{s}^{*}(\lambda)}{\mathrm{i}\lambda}
+a_{2}\frac{\widetilde{v}_{d}^{*}(\lambda)}{(\mathrm{i}\lambda)^{2}}-\left<y,v\right>(1+\frac{\mathrm{i}\alpha}{\lambda^2}m(\lambda))=0,\label{yv}
\end{equation*}
where $\widetilde{v}_{c}^{*}(\lambda)$, $\widetilde{v}_{s}^{*}(\lambda)$, $\widetilde{v}_{c}^{*}(\lambda)$, $m(\lambda)$ are defined in (\ref{vstar}) and (\ref{mlambda}).
This together with the periodic boundary conditions (\ref{bc}) implies that $z=\lambda^3$ is an eigenvalue of the operator $L_{\alpha}$ if and only if the linear equations system for $a_{0}$, $a_{1}$, $a_{2}$ and $\left<y,v\right>$
\begin{equation*}
M(\alpha, \lambda)
\left(\begin{array}{c}
a_{0}\\
a_{1}\\
a_{2}\\
\left<y,v\right>
\end{array}
\right)=\left(\begin{array}{c}
0\\
0\\
0\\
0
\end{array}
\right)\label{123yv}
\end{equation*}
has non-trivial solutions. Hence, $\Delta(\alpha,\lambda)$ is the characteristic function of the operator $L_{\alpha}$. On account of Lemma \ref{wcsd} and Lemma \ref{mmstar}, a simple manipulation gives
\begin{eqnarray*}
\Delta(\alpha,\lambda)\hh &=&\hh (1+\frac{\mathrm{i}\alpha}{\lambda^2}m(\lambda))\Delta(0,\lambda)-\frac{\mathrm{i}\alpha}{\lambda^2}\left\{\int_{0}^{1}d(\mathrm{i}\lambda(1-t))v(t)\mathrm{d}t
\left[\widetilde{v}_{c}^{*}(\lambda)\left(c(-\mathrm{i}\lambda)\right.\right.\right.\\
\hh &&\hh \left.\left.-2c(\mathrm{i}\lambda)+1\right)+\widetilde{v}_{s}^{*}(\lambda)(d(\mathrm{i}\lambda)+d(-\mathrm{i}\lambda))+\widetilde{v}_{d}^{*}(\lambda)(s(-\mathrm{i}\lambda)+s(\mathrm{i}\lambda))\right]\\
\hh &&\hh +\int_{0}^{1}s(\mathrm{i}\lambda(1-t))v(t)\mathrm{d}t\left[\widetilde{v}_{c}^{*}(\lambda)(s(-\mathrm{i}\lambda)+s(\mathrm{i}\lambda))
+\widetilde{v}_{s}^{*}(\lambda)\left(c(-\mathrm{i}\lambda)-2c(\mathrm{i}\lambda)\right.\right.\\
\hh &&\hh \left.\left.+1\right)+\widetilde{v}_{d}^{*}(\lambda)(d(\mathrm{i}\lambda)+d(-\mathrm{i}\lambda))\right]
+\int_{0}^{1}c(\mathrm{i}\lambda(1-t))v(t)\mathrm{d}t\left[\widetilde{v}_{c}^{*}(\lambda)\left(d(\mathrm{i}\lambda)\right.\right.\\
\hh &&\hh \left.\left.\left.+d(-\mathrm{i}\lambda)\right)+\widetilde{v}_{s}^{*}(\lambda)(s(-\mathrm{i}\lambda)+s(\mathrm{i}\lambda))
+\widetilde{v}_{d}^{*}(\lambda)(c(-\mathrm{i}\lambda)-2c(\mathrm{i}\lambda)+1)\right]\right\}\\
\hh &=&\hh \Delta(0,\lambda)+\frac{\mathrm{i}\alpha}{\lambda^2}\left\{m(\lambda)[3c(\mathrm{i}\lambda)-1]+m^{*}(\lambda)[3c(-\mathrm{i}\lambda)-1]\right.\\
\hh &&\hh \left.-\widetilde{v}_{c}^{*}(\lambda)\widetilde{w}_{d}(-\lambda)-\widetilde{v}_{c}(\lambda)\widetilde{w}_{d}^{*}(-\lambda)-\widetilde{v}_{s}^{*}(\lambda)\widetilde{w}_{s}(-\lambda)\right.\\
\hh &&\hh \left.-\widetilde{v}_{s}(\lambda)\widetilde{w}_{s}^{*}(-\lambda)-\widetilde{v}_{d}^{*}(\lambda)\widetilde{w}_{c}(-\lambda)-\widetilde{v}_{d}(\lambda)\widetilde{w}_{c}^{*}(-\lambda)\right\}.
\end{eqnarray*}
Using the definition of $F(\lambda)$ in (\ref{flambda}), we obtain (\ref{deltaa}). Furthermore, due to (iii) in Lemma \ref{entire} and (\ref{d0pp}), the identities in (\ref{dapp}) can be proved.
\qed

Now we explore the spectrum $\sigma(L_{\alpha})$ of the operator $L_{\alpha}$ via the expression of the resolvent $R_{L_{\alpha}}(z)$ (\ref{rla}).
\begin{lem}
\label{vpneq0}
For $v\in L^2_{\mathbb{C}}(0,1)$, $n\in\mathbb{Z}$, if $v_{n}=\left<v,u_{n}\right>\neq 0$, then $z_{n}\in \sigma (L_{0})$ is not the eigenvalue of the operator $L_{\alpha}$.
\end{lem}
\Proof
For $z_{n}\in \sigma (L_{0})$, $n\in\mathbb{Z}$, if $v_{n}\neq 0$, we deduce that the residue of the resolvent $R_{L_{\alpha}}(z)$ at $z_{n}$ is
\begin{equation*}
\text{Res}_{z_{n}}R_{L_{\alpha}}(z)f=\lim_{z\rightarrow z_{n}}( z_{n}-z)R_{L_{\alpha}}(z)f=f_{n}u_{n}-\frac{\alpha f_{n}\overline{v}_{n}}{\alpha v_{n}\overline{v}_{n}}\cdot {v}_{n}{u}_{n}=0,
\end{equation*}
hence $z_{n}\notin \sigma (L_{\alpha})$.\qed

Naturally, we can divide the set $\sigma (L_{0})$ into two disjoint subsets $\sigma (L_{0})=\sigma_{0}\cup\sigma_{1}$, where
\begin{equation}
\sigma_{0}:=\{z_{n}^{0}=z_{n}\in \sigma (L_{0})\vert v_{n}=0\},\quad\sigma_{1}:=\{z_{n}^{1}=z_{n}\in \sigma (L_{0})\vert v_{n}\neq0\}.\label{s0s1}
\end{equation}

\begin{prop}
\label{spectrumla}
For $v\in L^{2}_{\mathbb{C}}(0,1)$, the spectrum of the operator $L_{\alpha}$ is
\begin{equation*}
\sigma (L_{\alpha})=\sigma_{0}\cup\sigma_{2},
\end{equation*}
where
\begin{equation}
\sigma_{2}:=\{\mu_{n}\vert Q(\mu_{n})=0, n\in\mathbb{Z}\},\quad Q(z):=1+\alpha\sum\limits_{k\in\mathbb{Z}}\frac{\vert v_{k}\vert^{2}}{z_{k}-z}.\label{qz}
\end{equation}
Besides, the zeros of $Q(z)$ are real, simple and alternate with numbers $z_{n}^{1}\in \sigma_{1}$.
\end{prop}
\Proof
Lemma \ref{vpneq0} shows the resolvent $R_{L_{\alpha}}(z)$ does not have singularities at $z_{n}^{1}\in \sigma_{1}$, and $\sigma_{1}\nsubseteq \sigma (L_{\alpha})$. Owing to the division of the spectrum $\sigma (L_{0})$, we can simplify (\ref{rla}) as follows.
\begin{eqnarray}
R_{L_{\alpha}}(z)f\hh &=&\hh\sum_{n\in\mathbb{Z}}\frac{f_{n}}{z_{n}-z}u_{n}-
\alpha\frac{ \sum\limits_{k\in\mathbb{Z}}\frac{f_{k}\overline{v}_{k}}{z_{k}-z}}{1+\alpha \sum\limits_{k\in\mathbb{Z}}\frac{\vert v_{k}\vert^{2}}{z_{k}-z}}\cdot \sum\limits_{n\in\mathbb{Z}}\frac{v_{n}}{z_{n}-z}u_{n}\notag\\
\hh &=&\hh \sum\limits_{z_{n}^{0}\in \sigma_{0}}\frac{f_{n}}{z_{n}^{0}-z}u_{n}+\sum\limits_{z_{n}^{1}\in \sigma_{1}}\frac{f_{n}}{z_{n}^{1}-z}u_{n}
-\alpha\frac{ \sum\limits_{z_{n}^{1}\in \sigma_{1}}\frac{f_{n}\overline{v}_{n}}{z_{n}^{1}-z}}{Q(z)}\cdot \sum\limits_{z_{n}^{1}\in \sigma_{1}}\frac{v_{n}}{z_{n}^{1}-z}u_{n}\notag\\
\hh &=&\hh \sum\limits_{z_{n}^{0}\in \sigma_{0}}\frac{f_{n}}{z_{n}^{0}-z}u_{n}+\frac{1}{Q(z)}\sum\limits_{z_{n}^{1}\in \sigma_{1}}\frac{\left(f_{n}Q(z)-\alpha v_{n} \sum\limits_{z_{n}^{1}\in \sigma_{1}}\frac{f_{n}\overline{v}_{n}}{z_{n}^{1}-z}\right)u_{n}}{z_{n}^{1}-z}.\label{rela}
\end{eqnarray}
Thus, the zeros of $Q(z)$ and $z_{n}^{0}$ are the poles of $R_{L_{\alpha}}(z)$, and $\sigma (L_{\alpha})=\sigma_{0}\cup\sigma_{2}$. Let
\begin{equation}
G(z)=\sum\limits_{z_{n}^{1}\in \sigma_{1}}\frac{\vert v_{n}\vert^{2}}{z_{n}^{1}-z},\label{gz}
\end{equation}
then $Q(z)=1+\alpha G(z)$. For $z\in \mathbb{R}/\sigma_{1}$, the inequality $G^{\prime}(z)>0$ illustrates the function $G(z)$ is monotonically increasing, and $Q(z)=0$ has at most one real root in each interval $(z_{n-1}^{1},z_{n}^{1})$. On the other hand, the inequality $\lim\limits_{z\rightarrow z_{n-1}^1+}Q(z)\cdot \lim\limits_{z\rightarrow z_{n}^1-}Q(z)<0$ shows that $Q(z)$ has at least one real zero in each interval $(z_{n-1}^{1},z_{n}^{1})$. Therefore, the zeros of $Q(z)$ are real, simple and alternate with numbers $z_{n}^{1}\in \sigma_{1}$.
\qed

\begin{thm}
\label{spectrum}
For $v\in L^{2}_{\mathbb{C}}(0,1)$, the spectrum of the operator $L_{\alpha}$ is discrete and the multiplicity of each eigenvalue is not higher than $2$, besides, the number of eigenvalues of multiplicity $2$ are finite.
\end{thm}
\Proof
From Proposition \ref{spectrumla}, the spectrum of the operator $L_{\alpha}$ is simple except at the points in $\sigma_{0}\cap\sigma_{2}$, where its multiplicity equals $2$. We claim $\sigma_{0}\cap\sigma_{2}$ is a finite set. In fact, if $\sigma_{0}\cap\sigma_{2}=\{\hat{\mu}_{n}\}_{1}^{\infty}$, then according to $z_{n}=(2n\pi)^3$, we find $ \lim\limits_{n\rightarrow \infty }\hat{\mu}_{n}=\infty$ and the series in
\begin{equation}
Q(\hat{\mu}_{n})=1+\alpha\sum\limits_{z_{n}^{1}\in \sigma_{1}}\frac{\vert v_{n}\vert^{2}}{z_{n}^{1}-\hat{\mu}_n}\label{md}
\end{equation}
converges uniformly. Due to $Q(\hat{\mu}_{n})=0$ and passing to the limit as $n\rightarrow \infty$ in (\ref{md}), one obtains $1=0$. This leads to a contradiction.\qed

Now we deduce the expressions of the eigenfunctions.
\begin{prop}
\rm{(i)} For any $z_{n}\in \sigma_{0}$, the corresponding eigenfunction is
\begin{eqnarray*}
u_n(\alpha,x)=\mathrm{e}^{2n\pi\mathrm{i}x}.
\end{eqnarray*}
%\begin{eqnarray*}
%u(\alpha,x,\lambda)&=&-\frac{a_{0}\mathrm{i}}{a(\lambda)}\left\{c(\mathrm{i}\lambda(x-1))\widetilde{v}_{d}(\lambda)+s(\mathrm{i}\lambda(x-1))\widetilde{v}_{s}(\lambda)+d(\mathrm{i}\lambda(x-1))\widetilde{v}_{c}(\lambda)\right.\\
%&&+c(\mathrm{i}\lambda x)\widetilde{w}_{d}(-\lambda)+s(\mathrm{i}\lambda x)\widetilde{w}_{s}(-\lambda)+d(\mathrm{i}\lambda x)\widetilde{w}_{c}(-\lambda)\\
%&&\left.+[1-3c(\mathrm{i}\lambda)]\int_{0}^{x}d(\mathrm{i}\lambda(x-t))v(t)\mathrm{d}t+[1-3c(-\mathrm{i}\lambda)]\int_{x}^{1}d(\mathrm{i}\lambda(x-t))v(t)\mathrm{d}t\right\},
%\end{eqnarray*}
%where
%\begin{eqnarray*}
%a(\lambda)&=&-\int_{0}^{1}\left\{c(\mathrm{i}\lambda(1-t))[d(\mathrm{i}\lambda)+d(-\mathrm{i}\lambda)]+s(\mathrm{i}\lambda(1-t))[s(\mathrm{i}\lambda)+s(-\mathrm{i}\lambda)]\right.\\
%&&\left.+d(\mathrm{i}\lambda(1-t))[c(-\mathrm{i}\lambda)-2c(\mathrm{i}\lambda)+1]\right\}v(t)\mathrm{d}t.
%\end{eqnarray*}
\rm{(ii)} For each $\mu_{p}\in \sigma_{2}$, the corresponding eigenfunction is
\begin{equation*}
\widetilde{u}_{p}(\alpha,x)=\frac{1}{\sqrt {G^{\prime}(\mu_{p})}}\sum_{z_{n}^{1}\in \sigma_{1}}\frac{v_{n}}{z_{n}^{1}-\mu_{p}}u_{n}.
\end{equation*}
\end{prop}
\Proof
(i) For $z_{n}\in \sigma_{0}$, due to $v_{n}=\left<v,u_{n}\right>=0$, it yields $u_{n}$ satisfies the equation
\begin{equation*}
L_{\alpha}u_n=L_{0}u_{n}+\alpha \left<u_{n},v\right>v=z_{n}u_n
\end{equation*}
and the periodic boundary conditions (\ref{bc}). Thus, $u_n(\alpha,x)=u_n$.\\
(ii) The equality (\ref{rela}) shows the residue of $R_{L_{\alpha}}(z)f$ at $\mu_{p}\in\sigma_{2}$ is
\begin{eqnarray*}
\text{Res}_{\mu_{p}}R_{L_{\alpha}}(z)f= \lim_{z\rightarrow \mu_{p}}( \mu_{p}-z)R_{L_{\alpha}}(z)f= \frac{\alpha}{Q^{\prime}(\mu_{p})}\sum_{z_{n}^{1}\in \sigma_{1}}\frac{v_{n}}{z_{n}^{1}-\mu_{p}}(\sum_{z_{k}^{1}\in \sigma_{1}}\frac{f_{k}\overline{v}_{k}}{z_{k}^{1}-\mu_{p}})u_{n}.
\end{eqnarray*}
From the definitions of $Q(z)$ (\ref{qz}) and $G(z)$ (\ref{gz}), we see
\begin{equation*}
Q^{\prime}(\mu_{p})=\alpha G^{\prime}(\mu_{p}), \quad G^{\prime}(\mu_{p})>0.
\end{equation*}
Denote
\begin{equation*}
\widetilde{u}_{p}=\frac{1}{\sqrt {G^{\prime}(\mu_{p})}}\sum_{z_{n}^{1}\in \sigma_{1}}\frac{v_{n}}{z_{n}^{1}-\mu_{p}}u_{n},
\end{equation*}
then a simple calculation yields $\Vert \widetilde{u}_{p} \Vert_{L^{2}}=1$ and $\widetilde{u}_{p}$ satisfies
\begin{eqnarray*}
L_{\alpha}\widetilde{u}_{p}\hh &=&\hh L_{0}\widetilde{u}_{p}+\alpha\left<\widetilde{u}_{p},v\right>v\\
\hh &=&\hh \frac{1}{\sqrt {G^{\prime}(\mu_{p})}}\sum_{z_{n}^{1}\in \sigma_{1}}\frac{v_{n}}{z_{n}^{1}-\mu_{p}}z_{n}^{1}u_{n}+\frac{\alpha}{\sqrt {G^{\prime}(\mu_{p})}}\sum\limits_{z_{n}^{1}\in \sigma_{1}}\frac{\vert v_{n}\vert^{2}}{z_{n}^{1}-\mu_{p}}v\\
\hh &=&\hh \frac{1}{\sqrt {G^{\prime}(\mu_{p})}}\sum_{z_{n}^{1}\in \sigma_{1}}\frac{v_{n}}{z_{n}^{1}-\mu_{p}}(z_{n}^{1}-\mu_{p}+\mu_{p})u_{n}+\frac{\alpha}{\sqrt {G^{\prime}(\mu_{p})}}\sum\limits_{z_{n}^{1}\in \sigma_{1}}\frac{\vert v_{n}\vert^{2}}{z_{n}^{1}-\mu_{p}}v\\
\hh &=&\hh \frac{1}{\sqrt {G^{\prime}(\mu_{p})}}\sum\limits_{z_{n}^{1}\in \sigma_{1}}(v_{n}+\frac{v_{n}\mu_{p}}{z_{n}^{1}-\mu_{p}})u_{n}+\frac{\alpha}{\sqrt {G^{\prime}(\mu_{p})}}\sum\limits_{z_{n}^{1}\in \sigma_{1}}\frac{\vert v_{n}\vert^{2}}{z_{n}^{1}-\mu_{p}}v\\
\hh &=&\hh \frac{1}{\sqrt {G^{\prime}(\mu_{p})}}\left[\sum\limits_{z_{n}\in \sigma(L_{0})}v_{n}u_{n}+{\alpha}\sum\limits_{z_{n}^{1}\in \sigma_{1}}\frac{\vert v_{n}\vert^{2}}{z_{n}^{1}-\mu_{p}}v\right]+\frac{\mu_{p}}{\sqrt {G^{\prime}(\mu_{p})}}\sum_{z_{n}^{1}\in \sigma_{1}}\frac{v_{n}}{z_{n}^{1}-\mu_{p}}u_{n}\\
\hh &=&\hh \frac{1}{\sqrt {G^{\prime}(\mu_{p})}}\left[1+{\alpha}\sum\limits_{z_{n}^{1}\in \sigma_{1}}\frac{\vert v_{n}\vert^{2}}{z_{n}^{1}-\mu_{p}}\right]v+{\mu_{p}}\widetilde{u}_{p}\\
\hh &=&\hh {\mu_{p}}\widetilde{u}_{p}
\end{eqnarray*}
and the periodic conditions (\ref{bc}). Therefore, $\widetilde{u}_{p}(\alpha,x)=\widetilde{u}_{p}$ is an eigenfunction of $L_{\alpha}$ with respect to $\mu_{p}$.
\qed
\begin{rem}
For $z_n\in\sigma_0\cap \sigma_2$, the eigenfunctions $u_{n}(\alpha,x)$ and $\widetilde{u}_{n}(\alpha,x)$ are linearly independent.
\end{rem}

\section{The inverse problem}
In this section, we consider the inverse problems for the operator $L_{\alpha}$. Firstly, we give the Ambarzumyan-type theorem as follows.
\begin{thm}
\label{ath}
For $v\in L_{\mathbb{C}}^2(0,1)$, if $\sigma(L_\alpha)=\sigma(L_0)$, then $v=0$.
\end{thm}
\Proof
Since the function $v\in L^{2}_{\mathbb{C}}(0,1)$ can be expressed by its Fourier series, one has
\begin{equation*}
v(x)=\sum_{n\in\mathbb{Z}}v_{n}\mathrm{e}^{2n\pi \mathrm{i}x}.
\end{equation*}
So the potential $v(x)$ can be reconstructed if the Fourier coefficients $\{v_n\}$ are given. From Proposition \ref{spectrumla}, we find $\sigma(L_\alpha)=\sigma_0\cup\sigma_2$. It follows from $\sigma(L_\alpha)=\sigma(L_0)$ that $\sigma(L_\alpha)=\sigma_0=\sigma(L_0)$. Then, owing to the definition of $\sigma_0$ in (\ref{s0s1}), one has $v_n=0$ for $n\in\mathbb{Z}$. Therefore, the Fourier series $v(x)=\sum\limits_{n\in\mathbb{Z}}v_{n}\mathrm{e}^{2n\pi \mathrm{i}x}$ implies $v=0$.
\qed

In the following, we investigate the spectral data for recovering the non-local potential $\{\alpha,v\}$.
\begin{lem}
The function $Q(z)$ $(\ref{qz})$ satisfies the following identity
\begin{equation}
Q(\lambda^{3})=\frac{\Delta(\alpha, \lambda)}{\Delta(0, \lambda)}.\label{qzd0da}
\end{equation}
\end{lem}
\Proof
From Lemma \ref{rlaa} and (\ref{qz}), we conclude
\begin{equation*}
Q(\lambda^{3})=1+\alpha\left<R_{L_{0}}(\lambda^{3})v,v\right>,
\end{equation*}
then due to expression of $R_{L_{0}}(\lambda^{3})$ in Lemma \ref{rel0}, (\ref{deltaa}) and (\ref{flambda}), we have
\begin{eqnarray*}
\left<R_{L_{0}}(\lambda^{3})v,v\right>\hh &=&\hh \frac{-\mathrm{i}}{\lambda^{2}\Delta(0,\lambda)}\left\{\left<\int_{0}^{x}d(\mathrm{i}\lambda(x-t))[1-3c(\mathrm{i}\lambda)]v(t)\mathrm{d}t
\right.\right.\\
\hh &&\hh +\int_{0}^{1}\left[d(-\mathrm{i}\lambda t)c(\mathrm{i}\lambda(x-1))+s(-\mathrm{i}\lambda t)s(\mathrm{i}\lambda(x-1))+c(-\mathrm{i}\lambda t)d(\mathrm{i}\lambda(x-1))\right.\\
\hh &&\hh \left.+d(\mathrm{i}\lambda x)c(\mathrm{i}\lambda(1-t))+s(\mathrm{i}\lambda x)s(\mathrm{i}\lambda(1-t))+c(\mathrm{i}\lambda x)d(\mathrm{i}\lambda(1-t))\right]v(t)\mathrm{d}t\\
\hh &&\hh \left.\left.+\int_{x}^{1}d(\mathrm{i}\lambda(x-t))[1-3c(-\mathrm{i}\lambda)]v(t)\mathrm{d}t,v\right>\right\}\\
\hh &=&\hh \frac{-\mathrm{i}}{\lambda^{2}\Delta(0,\lambda)}\left[\widetilde{w}^{*}_{c}(-\lambda)\widetilde{v}_{d}(\lambda)+\widetilde{w}^{*}_{s}(-\lambda)\widetilde{v}_{s}(\lambda)
+\widetilde{w}^{*}_{d}(-\lambda)\widetilde{v}_{c}(\lambda)+(1-3c(\mathrm{i}\lambda))m(\lambda)\right.\\
\hh &&\hh \left.\widetilde{v}^{*}_{c}(\lambda)\widetilde{w}_{d}(-\lambda)+\widetilde{v}^{*}_{s}(\lambda)\widetilde{w}_{s}(-\lambda)+\widetilde{v}^{*}_{d}(\lambda)\widetilde{w}_{c}(-\lambda)
+(1-3c(-\mathrm{i}\lambda))m^{*}(\lambda)\right]\\
\hh &=&\hh \frac{\mathrm{i}}{\lambda^{2}\Delta(0,\lambda)}[F^{*}(\lambda)+F^{*}(\lambda)]\\
\hh &=&\hh \frac{\Delta(\alpha, \lambda)-\Delta(0, \lambda)}{\alpha\Delta(0, \lambda)}.
\end{eqnarray*}
Therefore,
\begin{equation*}
Q(\lambda^{3})=1+\alpha\frac{\Delta(\alpha, \lambda)-\Delta(0, \lambda)}{\alpha\Delta(0, \lambda)}=\frac{\Delta(\alpha, \lambda)}{\Delta(0, \lambda)}.
\end{equation*}
\qed
\begin{lem}
The multiplicative expansions of characteristic functions $\Delta(0, \lambda)$ and $\Delta(\alpha, \lambda)$ are
\begin{equation}
\Delta(0, \lambda)=-8\mathrm{i}\lambda^{3}\prod _{n=1}^{\infty}\left(1-\frac{\lambda^{6}}{(2n\pi)^{6}}\right)=6\sum_{n=1}^{\infty}\frac{\left(\mathrm{i}\lambda\right)^{6n-3}}{(6n-3)!},\label{delta0hd}
\end{equation}
\begin{equation}
\Delta(\alpha, \lambda)=\left\{
\begin{array}{ll}
 {c_{0}\prod \limits_{n\in\mathbb{Z}}\left(1-\frac{\lambda^{3}}{\lambda_{n}^{3}(\alpha)}\right)},&\text{for }0\notin\sigma(L_{\alpha}),\\
 {c_{1}\lambda^{3}\prod\limits_{n\in\mathbb{Z},\lambda_{n}(\alpha)\neq 0}\left(1-\frac{\lambda^{3}}{\lambda_{n}^{3}(\alpha)}\right)},&\text{for }0\in \sigma(L_{\alpha}),
\end{array}
\right.\label{deltaahd}
\end{equation}
respectively, where
\begin{eqnarray}
c_0\hh &=&\hh i\alpha\left\vert\int_{0}^{1}v(t)\mathrm{d}t\right\vert^{2},\label{c0}\\
c_1\hh &=&\hh -\mathrm{i}\left[1+\alpha\Im\left( \int_0^1\int_0^x(x-t)^2v(t)\overline{v}(x)\mathrm{d}t\mathrm{d}x+2\int_0^1\int_0^1v(x)\overline{v}(t)\frac{(1-t+x)^5}{5!}\mathrm{d}t\mathrm{d}x\right)\right],\label{c1}
\end{eqnarray}
$\lambda_{n}(\alpha)$ represents the zero of $\Delta(\alpha, \lambda)$ and the notation $\Im \beta$ denotes the imaginary part of $\beta$.
\end{lem}
\Proof
Due to (\ref{Delta0}), we find
\begin{equation*}
\Delta(0, \lambda)=-8\mathrm{i}\sin\frac{\lambda}{2}\sin\frac{\omega\lambda}{2}\sin\frac{\omega^{2}\lambda}{2}.
\end{equation*}
From the Hadamard theorem on factorization \cite{entire}, we know
\begin{equation*}
\sin z=z\prod\limits_{k=1}^{\infty}\left(1-\frac{z^{2}}{(n\pi)^{2}}\right),
\end{equation*}
and hence we can obtain the multiplicative expansions of characteristic functions $\Delta(0, \lambda)$ (\ref{delta0hd}).

According to (\ref{dapp}), we find that $\omega^{k}\lambda_{n}(\alpha)$, $k=1,2,3$ are the roots of $\Delta(\alpha, \lambda)=0$. Then, due to the Hadamard theorem on factorization \cite{entire}, suppose
\begin{equation}
\Delta(\alpha, \lambda)=c\lambda^{k}\mathrm{e}^{b\lambda}\prod\limits_{n\in\mathbb{Z},\lambda_{n}(\alpha)\neq 0}\left(1-\frac{\lambda^{3}}{\lambda_{n}^{3}(\alpha)}\right).\label{hadama}
\end{equation}
From the definition of $F(\lambda)$ in (\ref{flambda}) and the Taylor series of $c(z)$, $s(z)$ and $d(z)$, we find
\begin{eqnarray*}
&&\Delta(0, \lambda)=3[c(\mathrm{i}\lambda)+c(-\mathrm{i}\lambda)]=6\sum_{n=1}^{\infty}\frac{\left(\mathrm{i}\lambda\right)^{6n-3}}{(6n-3)!},\\
&&F(\lambda)=a_{1}(v)\lambda^{2}+a_{2}(v)\lambda^{5}+\cdots+a_n(v)\lambda^{3n-1}+\cdots,
\end{eqnarray*}
where
\begin{eqnarray*}
a_{1}(v)\hh &=&\hh -2\left<\int_0^x\frac{(x-t)^2}{2}v(t)\mathrm{d}t,v(x)\right>+\left<v(x),1\right>\left<\frac{t^2}{2},v(1-t)\right>\\
\hh &&\hh +\left<v(x),x\right>\left<t,v(1-t)\right>+\left<v(x),\frac{x^2}{2}\right>\left<1,v(1-t)\right>,
\end{eqnarray*}
\begin{eqnarray*}
a_{2}(v)\hh &=&\hh 3\cdot\frac{\mathrm{i}^3}{3!}\left<\int_0^x\frac{\mathrm{i}^2(x-t)^2}{2}v(t)\mathrm{d}t,v(x)\right>
-\left<v(x),1\right>\left<\frac{(-\mathrm{i}t)^5}{5!},v(1-t)\right>\\
\hh &&\hh -\left<v(x),\frac{(\mathrm{i}x)^3}{3!}\right>\left<\frac{(-\mathrm{i}t)^2}{2!},v(1-t)\right>
-\left<v(x),\frac{(\mathrm{i}x)^4}{4!}\right>\left<-\mathrm{i}t,v(1-t)\right>\\
\hh &&\hh -\left<v(x),\mathrm{i}x\right>\left<\frac{(-\mathrm{i}t)^4}{4!},v(1-t)\right>
-\left<v(x),\frac{(\mathrm{i}x)^2}{2!}\right>\left<\frac{(-\mathrm{i}t)^3}{3!},v(1-t)\right>\\
\hh &&\hh -\left<v(x),\frac{(\mathrm{i}x)^5}{5!}\right>\left<1,v(1-t)\right>,
\end{eqnarray*}
and $a_n(v)$ can be calculated similarly. Therefore, the equality (\ref{deltaa}) shows
\begin{eqnarray*}
\Delta(\alpha, \lambda)\hh &=&\hh \mathrm{i}\alpha[a_{1}(v)+\overline{a}_{1}(v)]+\left[\mathrm{i}\alpha(a_{2}(v)+\overline{a}_{2}(v))+\mathrm{i}^{3}\right]\lambda^{3}+\cdots+\left[\mathrm{i}\alpha\left(a_{2k+1}(v)\right.\right.\\
\hh&&\hh\left.\left.+\overline{a}_{2k+1}(v)\right)\right]\lambda^{6k-3}+\left[\mathrm{i}\alpha(a_{2k}(v)+\overline{a}_{2k}(v))+\left(\frac{6\mathrm{i}^{6k-3}}{(6k-3)!}\right)\right]\lambda^{6k}+\cdots\\
\hh &=&\hh c_0+c_1\lambda^{3}+\cdots+c_n\lambda^{3n}+\cdots,
\end{eqnarray*}
here, $c_0$, $c_1$ are defined in (\ref{c0}), (\ref{c1}), respectively.

If $z=0\notin \sigma(L_{\alpha})$, then $c_0\neq0$. According to (\ref{hadama}), we find $k=0$. Thus, due to $\Delta^{\prime}(\alpha, 0)=0$, we see $b=0$ and
\begin{equation*}
\Delta(\alpha, \lambda)=c_0\prod \limits_{n\in\mathbb{Z}}\left(1-\frac{\lambda^{3}}{\lambda_{n}^{3}(\alpha)}\right).
\end{equation*}
If $z=0\in \sigma(L_{\alpha})$, then $c_0=0$. The equality (\ref{hadama}) shows $k=3$. Then the identity $\Delta^{\prime}(\alpha, 0)=0$ yields $b=0$ and
\begin{equation*}
\Delta(\alpha, \lambda)=c_1\lambda^{3}\prod \limits_{n\in\mathbb{Z},\lambda_{n}(\alpha)\neq 0}\left(1-\frac{\lambda^{3}}{\lambda_{n}^{3}(\alpha)}\right).
\end{equation*}
\qed

\begin{lem}
\label{sd}
For $v\in L^{2}_{\mathbb{C}}(0,1)$, $\Vert v\Vert=1$, the spectra data $\{\alpha, \vert v_{p}\vert^{2}\}$ can be recovered by the spectra $\sigma(L_{0})$ and $\sigma(L_{\alpha})$.
\end{lem}
\Proof
(i) For $0\notin \sigma(L_{\alpha})$, we find $0\in \sigma_{1}$, then combining (\ref{qz}), (\ref{qzd0da}), (\ref{delta0hd}) and (\ref{deltaahd}), we find
\begin{eqnarray}
Q(z)\hh &=&\hh\frac{\Delta(\alpha,z^{\frac{1}{3}})}{\Delta(0,z^{\frac{1}{3}})}= 1+\alpha\sum\limits_{z_{n}^{1}\in\sigma_{1}}\frac{\vert v_{n}\vert^{2}}{z_{n}^{1}-z}\notag\\
\hh &=&\hh\frac{c_0 {\prod\limits_{n\in\mathbb{Z}}\left(1-\frac{z}{z_n(\alpha)}\right)}}{ {-8\mathrm{i}z\prod\limits_{n\in\mathbb{Z},z_n\neq0}\left(1-\frac{z}{z_n}\right)}}=\frac{ {c_0\prod\limits_{\mu_n\in\sigma_2}\left(1-\frac{z}{\mu_n}\right)}}{ {-8\mathrm{i}z\prod\limits_{z_n^1\in\sigma_1/\{0\}}\left(1-\frac{z}{z_n^1}\right)}},\label{qzz}
\end{eqnarray}
where $z_n(\alpha)=\lambda_{n}^3(\alpha)$. Hence,
\begin{equation*}
\frac{1}{c_0}=\lim_{y\rightarrow \infty}\frac{ {\prod\limits_{\mu_n\in\sigma_2}\left(1-\frac{\mathrm{i}y}{\mu_n}\right)}}{ {8y\prod\limits_{z_n^1\in\sigma_1/\{0\}}\left(1-\frac{\mathrm{i}y}{z_n^1}\right)}}.
\end{equation*}
Calculating the residue at the point $z_{n}^1$ in equality (\ref{qzz}), we have
\begin{eqnarray*}
&&\alpha \vert v_{n}\vert^{2}=-\frac{c_0\mathrm{i}}{8}\frac{z_n^1-\mu_{n}}{\mu_{n}}\prod\limits_{k\neq n}\frac{z_{k}^{1}}{\mu_{k}}\left(1-\frac{z_{k}^{1}-\mu_{k}}{z_{k}^{1}-z_{n}^1}\right),\quad n\neq 0,\\
&&\alpha \vert v_{0}\vert^2=\frac{c_0}{8\mathrm{i}}.
\end{eqnarray*}
Therefore, $\Vert v\Vert=1$ yields the spectra data $\{\alpha, \vert v_{p}\vert^{2}\}$.\\
(ii) For $0\in \sigma(L_{\alpha})$, we obtain $0\in \sigma_{0}$. It follows from (\ref{qz}), (\ref{qzd0da}), (\ref{delta0hd}) and (\ref{deltaahd}) that
\begin{eqnarray}
Q(z)\hh &=&\hh \frac{\Delta(\alpha,z^{\frac{1}{3}})}{\Delta(0,z^{\frac{1}{3}})}=1+\alpha\sum\limits_{z_{n}^{1}\in\sigma_{1}}\frac{\vert v_{n}\vert^{2}}{z_{n}^{1}-z}\notag\\
\hh &=&\hh \frac{c_1  {z\prod\limits_{n\in\mathbb{Z}, z_n(\alpha)\neq 0}\left(1-\frac{z}{z_n(\alpha)}\right)}}{ {-8\mathrm{i}z\prod\limits_{n\in\mathbb{Z},n\neq 0}\left(1-\frac{z}{z_n}\right)}}= \frac{ {c_1\prod\limits_{\mu_n\in\sigma_2}\left(1-\frac{z}{\mu_n}\right)}}{ {-8\mathrm{i}\prod\limits_{z_n^1\in\sigma_1}\left(1-\frac{z}{z_n^1}\right)}}.\label{qzz0}
\end{eqnarray}
Therefore,
\begin{equation*}
\frac{1}{c_1}=\lim_{y\rightarrow \infty}\frac{ {\prod\limits_{\mu_n\in\sigma_2}\left(1-\frac{\mathrm{i}y}{\mu_n}\right)}}{ {8\prod\limits_{z_n^1\in\sigma_1}\left(1-\frac{\mathrm{i}y}{z_n^1}\right)}}.
\end{equation*}
Calculating the residue at the point $z_{n}^1$ in equality (\ref{qzz0}), we obtain
\begin{equation*}
\alpha \vert v_{n}\vert^{2}=-\frac{c_1\mathrm{i}}{8}\frac{z_n^1-\mu_{n}}{\mu_{n}}\prod\limits_{k\neq n}\frac{z_{k}^{1}}{\mu_{k}}\left(1-\frac{z_{k}^{1}-\mu_{k}}{z_{k}^{1}-z_{n}^1}\right).
\end{equation*}
Then due to the fact $\Vert v\Vert=1$, the spectral data $\{\alpha, \vert v_{p}\vert^{2}\}$ is obtained.
\qed

\begin{thm}
\label{recover}
For $v\in L^{2}_{\mathbb{C}}(0,1)$, $\Vert v\Vert=1$, the non-local potential $\{\alpha, v\}$ can be recovered by four spectra $\sigma(L_{0})$, $\sigma(L_{\alpha}(v))$, $\sigma(L_{\alpha}(v+g))$ and $\sigma(L_{\alpha}(v+\mathrm{i}g))$, where $g(x)=1-x$.
\end{thm}
\Proof
%Since the function $v\in L^{2}_{\mathbb{C}}(0,1)$ can be expressed by its Fourier series, one has
%\begin{equation*}
%v(x)=\sum_{n\in\mathbb{Z}}v_{n}\mathrm{e}^{2n\pi \mathrm{i}x}.
%\end{equation*}
%So the potential $v(x)$ can be reconstructed if the Fourier coefficients $\{v_n\}$ are given.
The Lemma \ref{sd} shows the spectra data $\{\alpha, \vert v_{p}\vert^{2}\}$ can be recovered by $\sigma(L_{0})$ and $\sigma(L_{\alpha}(v))$. Similarly, we can calculate $\alpha\vert v_{n}+g_{n}\vert^2$ by $\sigma(L_{0})$ and $\sigma(L_{\alpha}(v+g))$. Let $\Re \beta$ denote the real part of $\beta$. Because of the fact that
\begin{equation*}
\alpha\vert v_{n}+g_{n}\vert^2=\alpha\left(\vert v_{n}\vert^2+2\Re (v_{n}\overline{g}_{n})+\vert g_{n}\vert^2\right),
\end{equation*}
and
\begin{equation*}
g_{n}=\left<g(x),\mathrm{e}^{2n\pi\mathrm{i}x}\right>=\left\{
\begin{array}{ll}
\frac{1}{2},&n=0,\\
-\frac{1}{2n\pi\mathrm{i}},&n\neq 0,
\end{array}
\right.
\end{equation*}
the numbers $\Re (v_{n}\overline{g}_{n})$ can be calculated by three spectra $\sigma(L_{0})$, $\sigma(L_{\alpha}(v))$, $\sigma(L_{\alpha}(v+g))$. Analogously, from $\sigma(L_{0})$, $\sigma(L_{\alpha}(v))$ and $\sigma(L_{\alpha}(v+\mathrm{i}g))$, we obtain $\Im (v_{n}\overline{g}_{n})$. Thus, $v_{n}\overline{g}_{n}$ and $v_{n}$ can be recovered by four spectra $\sigma(L_{0})$, $\sigma(L_{\alpha}(v))$, $\sigma(L_{\alpha}(v+g))$ and $\sigma(L_{\alpha}(v+\mathrm{i}g))$. Thereafter, the function $v(x)$ is recovered by its Fourier series $v(x)=\sum\limits_{n\in\mathbb{Z}}v_{n}\mathrm{e}^{2n\pi \mathrm{i}x}$.
\qed

\begin{thm}
\label{do}
For $v\in L^{2}_{\mathbb{C}}(0,1)$, $\overline{v}(1-x)=v(x)$, $\Vert v\Vert=1$, the non-local potential $\{\alpha, v\}$ can be recovered by three spectra $\sigma(L_{0})$, $\sigma(L_{\alpha}(v))$ and $\sigma(L_{\alpha}(v+h))$, where $h(x)=(x-\frac{1}{2})^{2}$.
\end{thm}
\Proof
From $\overline{v}(1-x)=v(x)$, we see
\begin{eqnarray*}
\overline{v}_{n}\hh&=&\hh\left<\mathrm{e}^{2n\pi\mathrm{i}x},v(x)\right>= \int_{0}^{1}\mathrm{e}^{2n\pi\mathrm{i}x}\overline{v}(x)\mathrm{d}x=\int_{0}^{1}\mathrm{e}^{2n\pi\mathrm{i}(1-t)}\overline{v}(1-t)\mathrm{d}t\\
\hh&=&\hh  \int_{0}^{1}\mathrm{e}^{-2n\pi\mathrm{i}x}v(x)\mathrm{d}x
= \left<v(x),\mathrm{e}^{2n\pi\mathrm{i}x}\right>= {v}_{n},
\end{eqnarray*}
so $v_{n}\in \mathbb{R}$. Then, following the same procedure as the proof of Theorem \ref{recover}, we find $\{\alpha, v_{n}^2\}$ can be recovered by $\sigma(L_{0})$ and $\sigma(L_{\alpha}(v))$, and $\alpha\vert v_{n}+h_{n}\vert^2=\alpha( v_{n}+h_{n})^2$ can be calculated by $\sigma(L_{0})$ and $\sigma(L_{\alpha}(v+h))$. Since
\begin{equation*}
h_{n}=\left<h(x),\mathrm{e}^{2n\pi\mathrm{i}x}\right>=\left\{
\begin{array}{ll}
\frac{1}{12},&n=0,\\
\frac{1}{2n^2\pi^2},&n\neq 0,
\end{array}
\right.
\end{equation*}
we find that the real numbers $v_{n}{h}_{n}$ and $v_{n}$ are unambiguously calculated by three spectra $\sigma(L_{0})$, $\sigma(L_{\alpha}(v))$, $\sigma(L_{\alpha}(v+h))$. Thus, the function $v(x)$ is reconstructed by its Fourier series $v(x)=\sum\limits_{n\in \mathbb{Z}}v_{n}\mathrm{e}^{2n\pi \mathrm{i}x}$.
\qed
\begin{cor}
\label{do1}
For $v\in L^{2}_{\mathbb{C}}(0,1)$, $\overline{v}(1-x)=-v(x)$, $\Vert v\Vert=1$, the non-local potential $\{\alpha, v\}$ can be recovered by three spectra $\sigma(L_{0})$, $\sigma(L_{\alpha}(v))$ and $\sigma(L_{\alpha}(v+\tilde {h}))$, where $\tilde {h}(x)=(x-\frac{1}{2})^{2}\mathrm{i}$.
\end{cor}
\Proof
The condition $\overline{v}(1-x)=-v(x)$ show that each Fourier coefficient $v_n$ is a pure imaginary number. Then, proceeding as in the proof of Theorem \ref{do}, we can prove this corollary.
\qed
\section*{Funding}
This research was supported by the Fundamental Research Funds for the Central Universities of Civil Aviation University of China (Grant No. 3122022061).

\end{document}